\documentclass[reqno]{jgokova}
\usepackage{epsfig,graphicx}

\newcommand{\ben}{\begin{enumerate}}
\newcommand{\een}{\end{enumerate}}
\newcommand{\be}{\begin{equation}}
\newcommand{\ee}{\end{equation}}
\newcommand{\bea}{\begin{eqnarray}}
\newcommand{\eea}{\end{eqnarray}}
\newcommand{\bc}{\begin{center}}
\newcommand{\ec}{\end{center}}

\newtheorem{thm}{Theorem}[section]

\theoremstyle{definition}

\theoremstyle{remark}

\newtheorem{Rm}{Remark}

\begin{document}

\setcounter{page}{1}

\title[Short title]{Homotopy $4$-spheres associated to an infinite order loose cork}
\author[]{Selman Akbulut}

\thanks{Partially supported by NSF grants DMS 0905917}

\address{ G\"{o}kova Geometry Topology Institute,  Mu\u{g}la, T\"{u}rk\.{i}ye}
 \email{akbulut.selman@gmail.com}

\begin{abstract}
We prove that the homotopy spheres $\Sigma_{n} = 
 -W\smile_{f^{n}}W$, formed by doubling the infinite order loose-cork $(W,f)$, by the iterates of the cork  automorphism $f: \partial W \to \partial W$, is $S^4$. To do this we first show that $\Sigma_{n} $ are obtained by Gluck twistings of $S^4$; then from this we show how to cancel $3$-handles of $\Sigma_{n}$ and identify it  by $S^{4}$.
\end{abstract}


\maketitle

\section{Introduction} \label{introduction} 

\begin{figure}[ht]  \begin{center}
 \includegraphics[width=.67\textwidth]{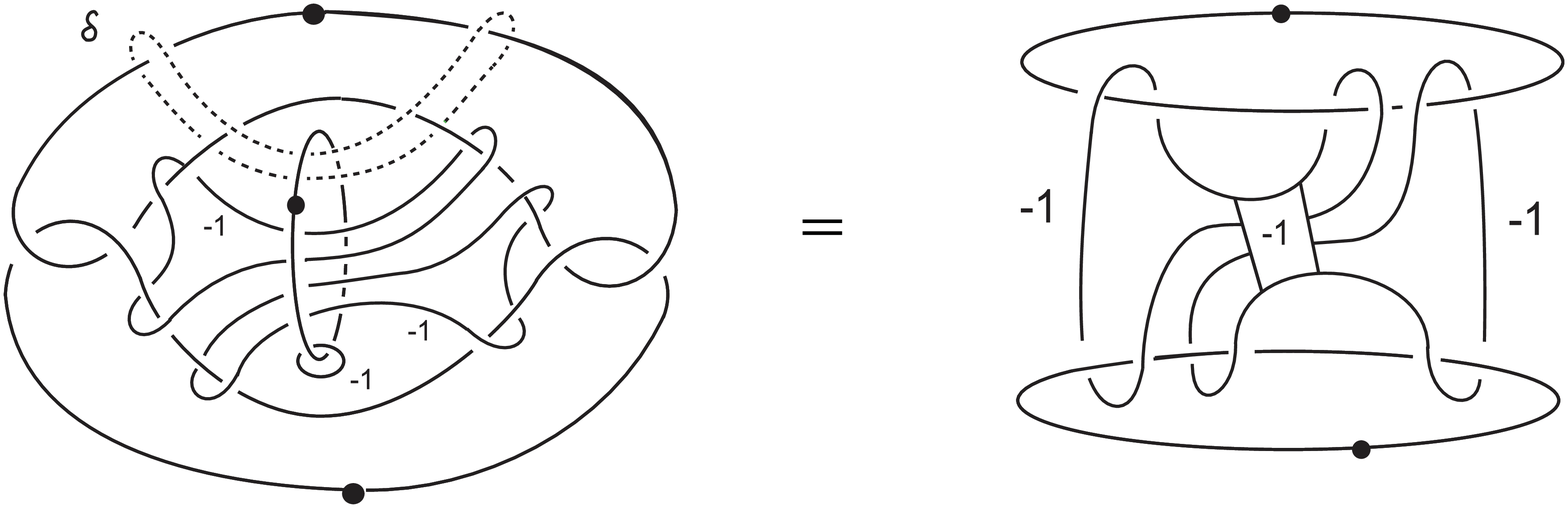}       
\caption{W}      \label{a1} 
\end{center}

\vspace{-.4in}

 \end{figure}

\begin{figure}[ht]  \begin{center}
 \includegraphics[width=.67\textwidth]{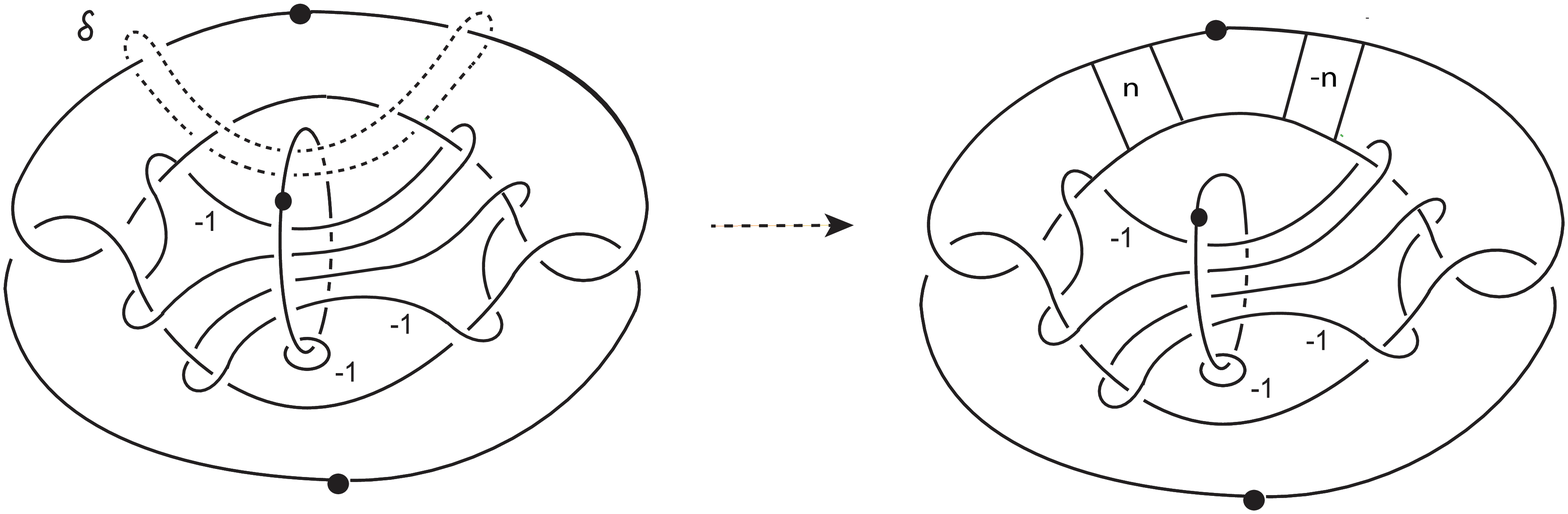}       
\caption{$\delta^{n} $ - move}      \label{a2} 
\end{center}
 \end{figure}

Let $(W,f)$ be the infinite order loose-cork of \cite{a1}, shown in Figure~\ref{a1}. As indicated in \cite{a1}, this $W$  can be identified with the one described in \cite{g1}. Recall that the diffeomorphism $f: \partial W \to \partial W$ here is given by the $\delta$-move along the curve $\delta$ of Figure~\ref{a2} as defined in \cite{a1}.  For simplicity we will refer the iterates $f^{n}$ of $f$ as $\delta^{n}$ or  $\delta $ -move.

 
  Even though the $\delta$-move diffeomorphism $f:\partial W \to \partial W$ by itself looks like an innocuous operation, when $W$ appears as a codimension zero submanifold $W\subset M^{4}$, the operation of cutting $W$ from $M$ and regluing with the composition map $f^{n}=f\circ f\circ...\circ f$ could result infinitely many different diffeomorphism types of $M$ (\cite{a1}). For example, reader can easily observe that if we attach a handle to $W$ along  $\gamma$, then the $\delta $-move operation applied to $W$ alters the position of   $\gamma$ as shown in Figure~\ref{a26}.  Consider the homotopy $4$-spheres  obtained by doubling of the contractible
manifold $W$ by the iterates $f^{n}=f\circ f\circ...\circ f$:
\begin{equation}\label{Sigma}\Sigma_{n}=- W\smile_{f^{n}} W \end{equation}
\noindent The obvious questions is whether $\Sigma_{n}$ are diffeomorphic to $S^4$,  or if this  family contains an exotic copy of $S^4$. The answer is given by:

\begin{thm} 
\label{Theorem1} \label{main}
Each $\Sigma_{n}$ is diffeomorphic to $S^4$.
\end{thm}

\vspace{.1in}

For  brevity call $\Sigma_{n}^{o}= \Sigma_{n}-B^4$.  Proof will proceed as follows: We will first show each $\Sigma_{n}$ is obtained by Gluck twisting $S^{4}$ along some knotted $S^2\subset S^{4}$, then find $3$-handle free handle pictures of $\Sigma_{n}$. We then cancel $3$-handles of $\Sigma_{n}$ by the trick, which was used in the solution of the Cappell-Shaneson homotopy sphere problem \cite{a2}. From this we will see that $\Sigma_{n}^{o}$  is obtained from the ribbon complement $Q=B^4-N(D^2)$  by attaching a $2$-handle along a knot $\gamma_{n}\subset \partial Q$, where $D\subset B^4$ is the standard ribbon bounded by $\partial D= K\# K$, where $K$ is the figure-$8$ knot, and $N(D)$ is the tubular neighborhood of $D$. 
\begin{equation} \label{withQ} \Sigma_{n}^{o} = Q\smile h_{\gamma_{n}}^{2}\end{equation}

Knots $\gamma_{n}$ are related each other by $f(\gamma_{n-1}) =\gamma_{n}$, where $f$ is a  $\delta$-move diffeomorphism $f:\partial Q \to \partial Q$. Next we show $\Sigma_{0}^{o} = B^4$, and 
$\Sigma_{n}^{0}$  is obtained by gluing $\Sigma_{n-1}^{0}$  to $ \partial \Sigma_{n-1}^{0}\times [0,1]$ along  $\partial \Sigma_{n-1}\times {0}$,  by a $\delta$-move diffeomorphism $g_{n-1}:  \partial \Sigma_{n-1}^{0} \to \partial \Sigma_{n-1}^{0}$, associated to $f$. 
$$\Sigma_{n}^{0} =\Sigma_{n-1}^{0}\smile_{g_{n-1}} \partial \Sigma_{n-1}^{0}\times [0,1] $$ Then this fact coupled with the fact that any $\delta$-move diffeomorphism  $g:S^3 \to S^3$ is isotopic to identity, finishes the proof by induction. 

\vspace{.1in}

\section{Construction}\label{construction}

Our first goal is to determine how the $\delta$-move diffeomorphism $f$ moves curves on the boundary $\partial W$ (see also \cite{a3}). This is important, because by using this we will construct the handlebody picture of the manifolds $\Sigma_{n}$,  by drawing the attaching circles of the dual handles of the upside down $-W$. This is a nontrivial task, because $\delta$-move is performed by first introducing and then canceling $2/3$-handle pairs. So the attaching $S^2$ of the $3$-handle might puncture the dual $2$-handle curves on $-\partial W$, forcing us to push them into the interior of $W$. To go around this problem, we will describe the $\delta$-move diffeomorphism in an alternative way, as a carving and uncarving operations, that is $1$- and $2$- handle exchanges in the interior (this is also referred as "dot and zero exchanges" in short).  This technique was exploited in \cite{a3}.

\vspace{.1in}

\newpage

In Figure~\ref{a3} we first replace dot with zero (turning $1$-handle to $2$-handle), then perform the  $2$-handle slide (indicated by the arrows), resulting the handlebody on the right. 
The reverse operation (i.e. going from right to left of the figure) can be obtained by first doing the $2$-handle slide, indicated by the dotted arrows, and then by replacing ``zero with dot". Here we also traced the dual circles to the $2$- handles during this operation (small red circles), where attaching $2$-handles to these circles gives the double of $W$, which we denoted by $\Sigma_{0}$. To construct the handlebody of $\Sigma_{n}$, we need to modify  $W \subset \Sigma_{n}$ along its boundary by a $\delta$-move.

 \vspace{.1in}
 
Figure~\ref{a4} indicates how the 
$\delta^{n}$-move $f: \partial W \to \partial W$ affects the dual handles (red) circles (figure drawn for $n=2$). Going back to the original Figure~\ref{a1} via the reverse $\delta  $-move (as indicated in  Figure~\ref{a3}) shows that the effect of the $\delta$-move on dual circles, is as  in  Figure~\ref{a5}.

\vspace{.1in}

Now comes a crucial point: A reader gazing at the first picture of Figure~\ref{a4} might conclude that $\delta$-move doesn't move the dual circles, because $n$ and $-n$ twists cancel each other. Here are two explanations: First of all, here we are dealing with circles-with-dots not framed circles, transferring twist across them has the affect of changing the carvings (i.e. changing the interiors). Secondly, the original $\delta$ move takes place on $\partial W$, not on the homotopy ball $\Sigma_{n}^{0}= W \smile [\mbox{dual} $ 2 \mbox{and} 3$-\mbox{handles}]$, that is $\delta$ may not be an on unknot on $\partial \Sigma_{n}^{0} $. Surprisingly, we can obtain Figure~\ref{a5} by performing $\delta$-move to $\Sigma_{n}^{0}$ (the first picture of Figure~\ref{a3}, with dual handles) by using the curve $d$ of Figure~\ref{a6}. This $d$ is in fact an unknot on 
$\partial (\Sigma_{n}^{0})$, which can be checked by the boundary correspondence of Figure~\ref{a3}. Also $d$  happens to be an unknot on $\partial W$, so we could use $d$ for the place of $\delta$, to serve for the dual purpose.

\vspace{.1in}

To sum up,  the first picture of Figure~\ref{a7} represents a handlebody of  $\Sigma_{n}$. Now it is easy to check that the middle dotted curve in the second picture of Figure~\ref{a7} is an unknot (to see this, do the reverse $\delta$-move go back to the first picture of Figure~\ref{a3}, and then observe that in the presence of the dual $2$-handles, the dotted circle becomes an unknot there). From this we see that $\Sigma_{n}$ is obtained from $S^4=\Sigma_{0}$ by Gluck twisting (this requires a simple check here, namely remove the dotted circle, and the $-1$ twist on the curves it links from the middle of Figure~\ref{a7}, then see that you get $S^{4}$). Now by using this unknot, we can attach a $2/3$-handle pair (the new $2$-handle is the $0$-framed dotted curve in the figure). Next we employ a trick , which was used solving the ``Cappell-Shaneson homotopy sphere problem'' (Figure 14.11 of \cite{a2}):  After the obvious  handle slide over the middle $0$-framed $2$-handle in Figure~\ref{a7}, we obtain the pictures of  Figure~\ref{a8},  where we can see two cancelling $1/2$ -handle pairs! The two $0$-framed middle $2$-handles cancel the two $1$-handles (represented with large dotted circles)! So this picture can be thought of a handlebody without $1$-handles, and hence turning it upside down we will get a handlebofy without $3$-handles! Having noted this, we can turn this handlebody upside down (as the process described in \cite{a2}). That is, we ignore the cancelled $1/2$ handle pairs, and carry the duals of the remaing $2$-handles to the boundary of $\# 3(S^{1}\times B^{3})$ by a diffeomorphism (duals to $2$-handles are indicated by the dashed little circles in Figure~\ref{a8}).

\vspace{.1in} 

Now our task is to find a diffeomorphism from the boundary of the pictures of Figure~\ref{a8}, to $\# 3(S^{1}\times B^{3})$ and carry the dual $2$-handles.  By applying the Figure~\ref{a3} boundary identification, we see that the boundary of the last picture of Figure~\ref{a8} can be identified with the boundary of the first picture of Figure~\ref{a9}, then the obvious isotopy gives the second picture of Figure~\ref{a9}. Note that we don't draw $3$- and $4$-handles here, the handldebodies of $\Sigma_{n}$ and $\Sigma_{n}^{0}$ will be drawn the same.  

\vspace{.1in} 

Again by applying the reverse boundary identification of Figure~\ref{a3} to Figure~\ref{a9} we get Figure~\ref{a10}, which is a handlebody picture of $\Sigma_{n}$, without $3$-handles! Finally the indicated simple handler slide gives Figure~\ref{a11} (the picture is drawn for $n=2$). To indicate how the pattern changes as we increse $n\to n+1$, in Figure~\ref{a12}, we drew $\Sigma_{n}$ for $n=1$. 

\vspace{.1in}  

 Now let us check the identification \ref{withQ} of Section~\ref{introduction}.   We will demonstrate a proof for $\Sigma_{1}$ (from this the reader can see the proof for the general case). For this we first isotope Figure~\ref{a12} to Figure~\ref{a13}, then do the handle slides and cancellations of the figures Figures~\ref{a13} 
$\leadsto $.. $\leadsto$ ~\ref{a18}, as indicated in the pictures. During these operations we trace the ribbon which the unknot $T$ of Figure~\ref{a13} bounds. In this figure this ribbon is the trivial ribbon bounding the unknot, where its ribbon move indicated with an unknotted arc in Figure~\ref{a15}. But during the handle slide Figure~\ref{a15} $\to$ Figure~\ref{a16} this trivial ribbon turns into the nontrivial ribbon $D$, mentioned above. By performing the ribbon move in Figure~\ref{a16}, along the indicated dotted arc, we get Figure~\ref{a17} (the dotted blue line of his figure is the dual of dotted red line of Figure~\ref{a16}). Then the  $1/2$ handle cancellation by using the $-1$ framed $2$-handle, gives Figure~\ref{a18} which is $\Sigma_{1}$. To see the general pattern, we can apply the same steps to  Figure~\ref{a11} rather than  Figure~\ref{a12}, then we see that we get Figure~\ref{a19} picture of $\Sigma_{2}$. Now the handlebody patterns of $\Sigma_{n}$ is as required in ~\ref{withQ} . 

\section{Rolling versus  carving}

  Notice that the loop $c\subset \partial M_{n}=S^3$, which links the ribbon in Figure~\ref{a18} (and in Figure~\ref{a19}) is the unknot in $S^3$. This is because doing $-1$ surgery to $c$ (which corresponds to putting $0$-framing on $c$ on the figure) gives $S^3$, hence by Property $P$ the loop $c$ must be the unknot. Now we can attach a cancelling $2/3$ handle pair to $\Sigma_{n}$  along $c$ (this corresponds to adding $+1$ framed $2$-handle to $c$).  This gives an alternative description of $\Sigma_{n}$ which contains a copy of $W$: 
  
 \begin{equation} \label{Wdecop} \Sigma_{n}^{0} = W\smile_{ f^{n}(\gamma)} h_{\gamma}^{2}\end{equation}
   
  This is because $W$ is in the form $W= Q \smile h^{2}_{c}$, where $Q=B^{4}-N(D)$ and $h_{c}$ is a $+1$-framed $2$ handle attached along $c$ (\cite{a1} Remark 1, and \cite{g1})), i.e. $\Sigma_{n}^{0}$  is obtained by attaching $2$-handle to $W$  along the $n$-th iterate of a loop $\gamma \subset \partial Q$ by some diffeomorphism $$f: \partial Q\to \partial Q$$  

\newpage

\begin{Rm} The handlebody picture of $\Sigma_{n}$  (Figures~\ref{a18} and \ref{a19}) shows that, by changing the carving, which $K\#K$ bounds in $B^4$ by a diffeomorphism will move the position of the $2$-handle $\gamma$ to the $2$-handle of Figure~\ref{a21}, which can easily be identified with $B^4$. This diffeomorphism is obtained by first moving the knot $K\# K$ by an isotopy $g_{t} :S^{3}\to S^{3}$ back to itself as indicated in Figure~\ref{a20} along the dotted arrow (i.e. {\it rolling} one of the factors of the connected sum over $K\#K$ back to itself), then letting $g_{1}(K\#K)$ bound the standard slice disk $D$ in $B^{4}$, which $K\#K$ bounds. Call this new ribbon disk $D'$. 
Recall that in \cite{a4} relatively exotic but diffeomorphic ribbon complements in $B^4$ were constructed. Here the  ribbon complements $D$ and $D'$ have the similar property (otherwise $W$ wouldn't be loose cork). \end{Rm}

\begin{Rm} 
Reader should compare this to the infinitely many absolutely exotic manifolds of \cite{a3}, which also decompose as \ref{Wdecop}. To study $\Sigma_{n}$ we have two options: (1) Either attach the rolled $2$-handle to the standard ribbon complement $B^{4}-D$ as in Figure~\ref{a19},  or (2) Attach the standard $2$ handle of Figure~\ref{a21} to the nonstandard ribbon complement $B^{4}-D'$, carved by rolling.  \end{Rm}

Note that the ``Dehn twist diffeomorphism $\partial W \to \partial W$ along an imbedded torus'',  discussed in \cite{g1} and \cite{rr}, corresponds to  $\delta$-move diffeomorphism along some $\delta \subset \partial W$. Patient reader can check this by tracing the steps outlined by Remark 1 of  \cite{a1}, one gets the identification of Figure~\ref{a24}. Then by doubling and connect summing the circle $\delta_1$ and the arc $\delta_2$ one can recover the position of $\delta$ on the right picture of $W$ in Figure~\ref{a24} (cf.\cite{a2}). This shows that the $\delta$ move of $W$ corresponds to Dehn twisting boundary of $W$ along the imbedded torus of Figure~\ref{a25}.  Also note that since any imbedded torus $T\subset S^3$ bounds a solid torus any Dehn twist diffeomorphism $S^3 \to S^3$ is isotopic to identity. Hence $\delta$ -move diffeomorphisms of $S^3$ are isotopic to identity.

\vspace{.1in}

\begin{figure}[ht]  \begin{center}
 \includegraphics[width=1\textwidth]{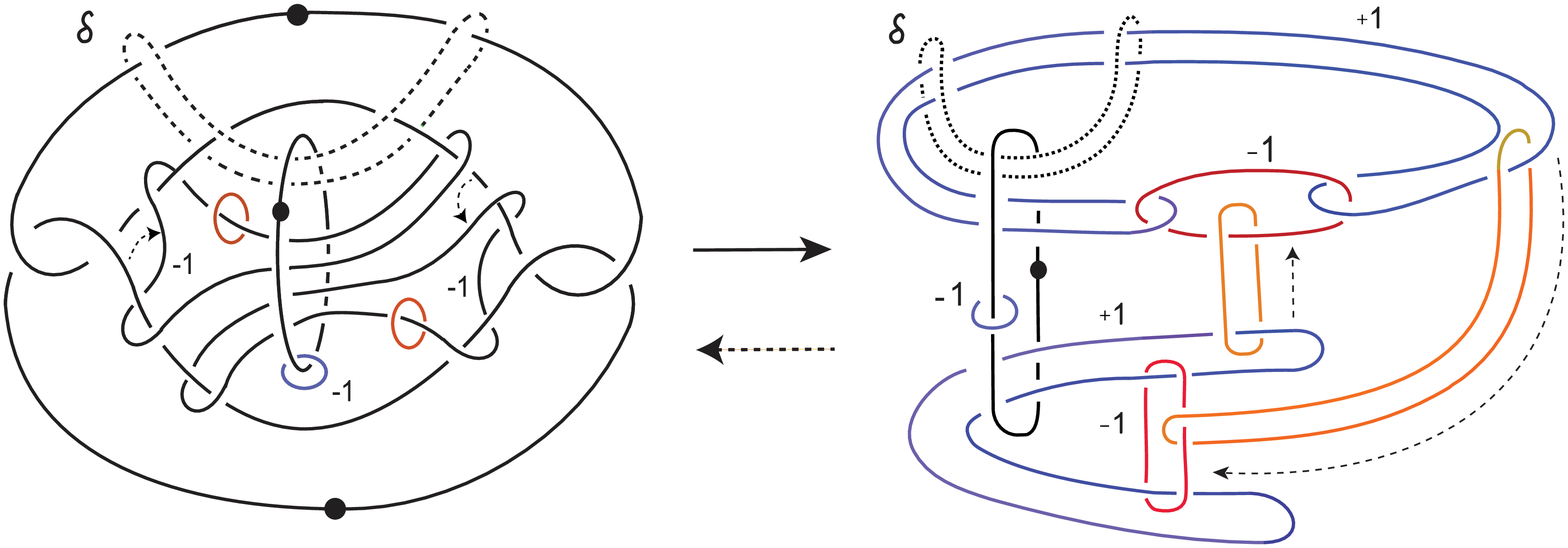}       
\caption{Changing the carvings}      \label{a3} 
\end{center}
 \end{figure}

\begin{figure}[ht]  \begin{center}
 \includegraphics[width=1\textwidth]{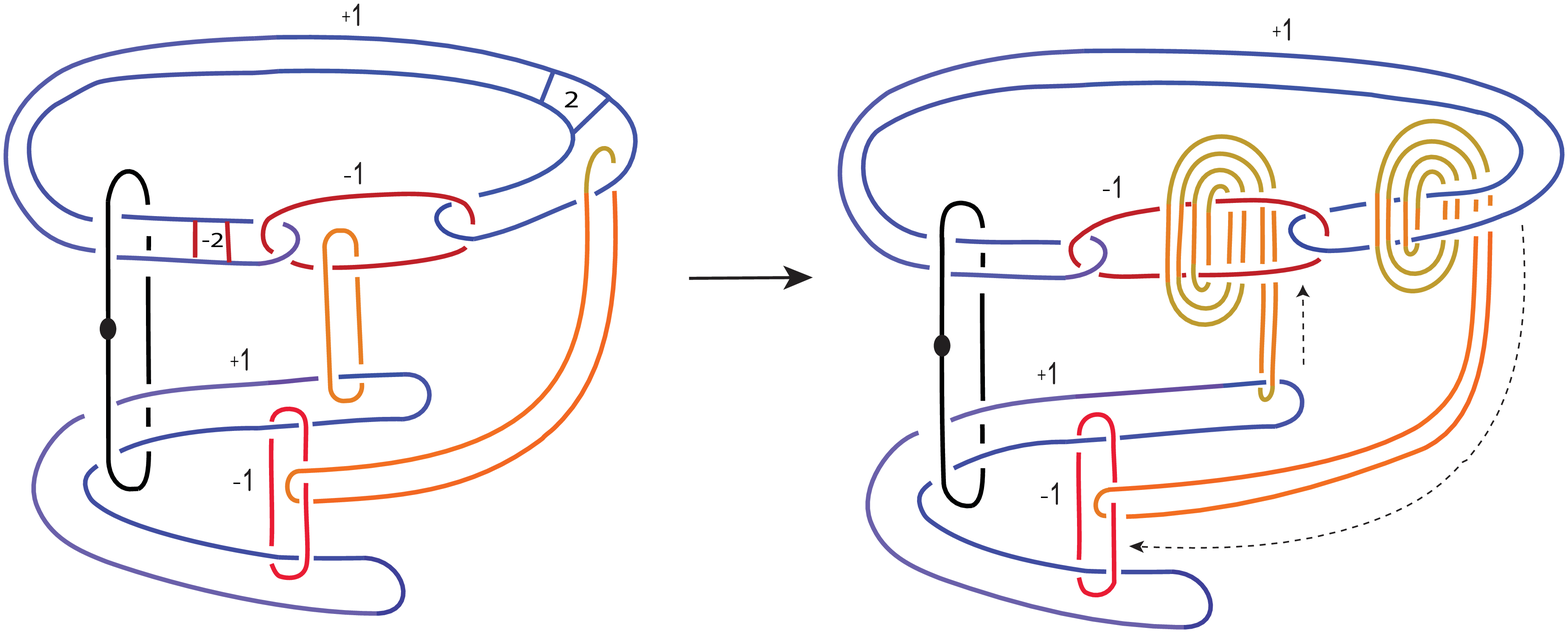}       
\caption{Affect of $\delta$-move on the boundary, n=2}   \label{a4} 
\end{center}
 \end{figure}

\begin{figure}[ht]  \begin{center}
 \includegraphics[width=.65\textwidth]{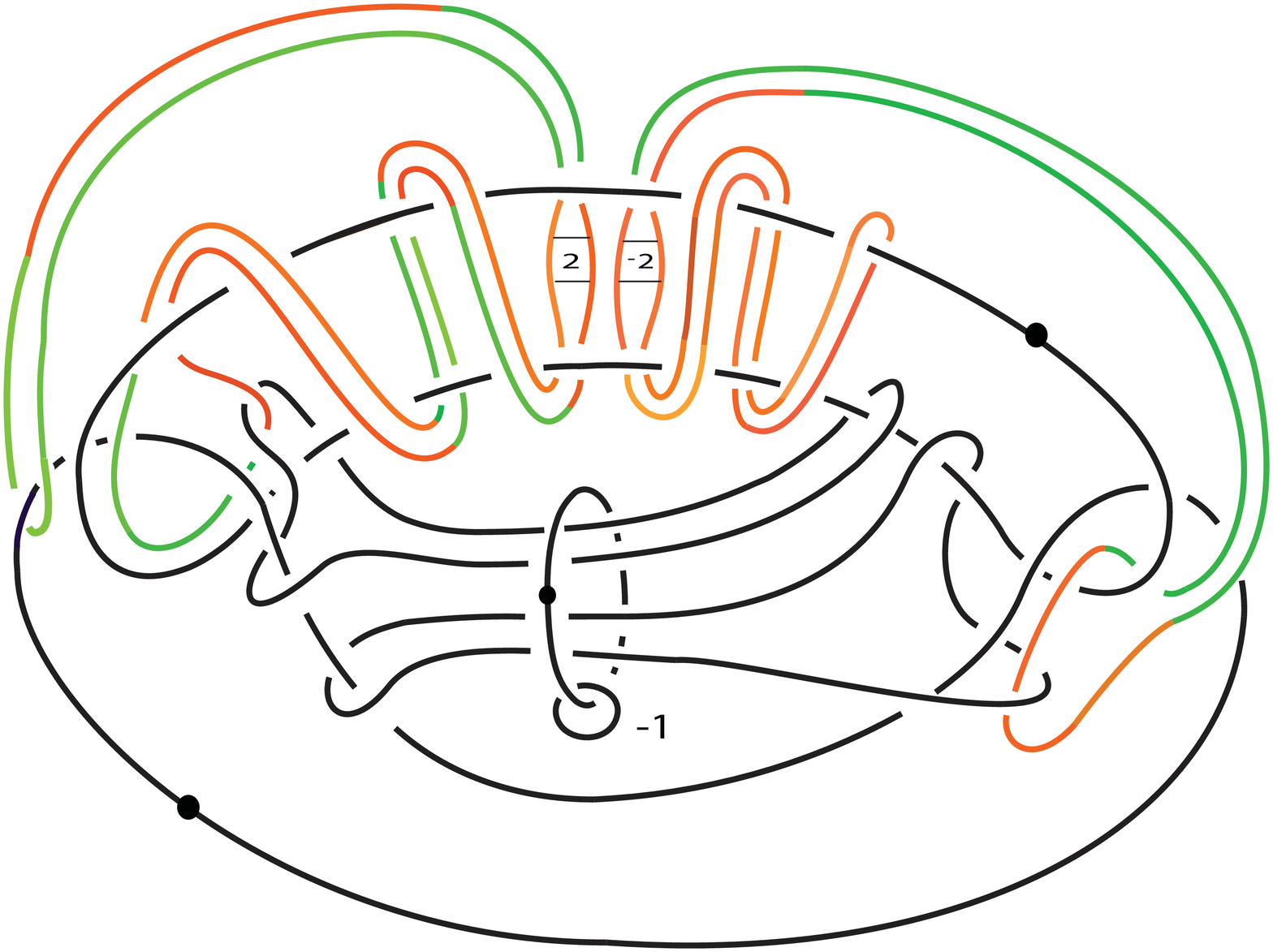}       
\caption{ $\Sigma_2$}   \label{a5} 
\end{center}
 \end{figure}

\begin{figure}[ht]  \begin{center}
 \includegraphics[width=.6\textwidth]{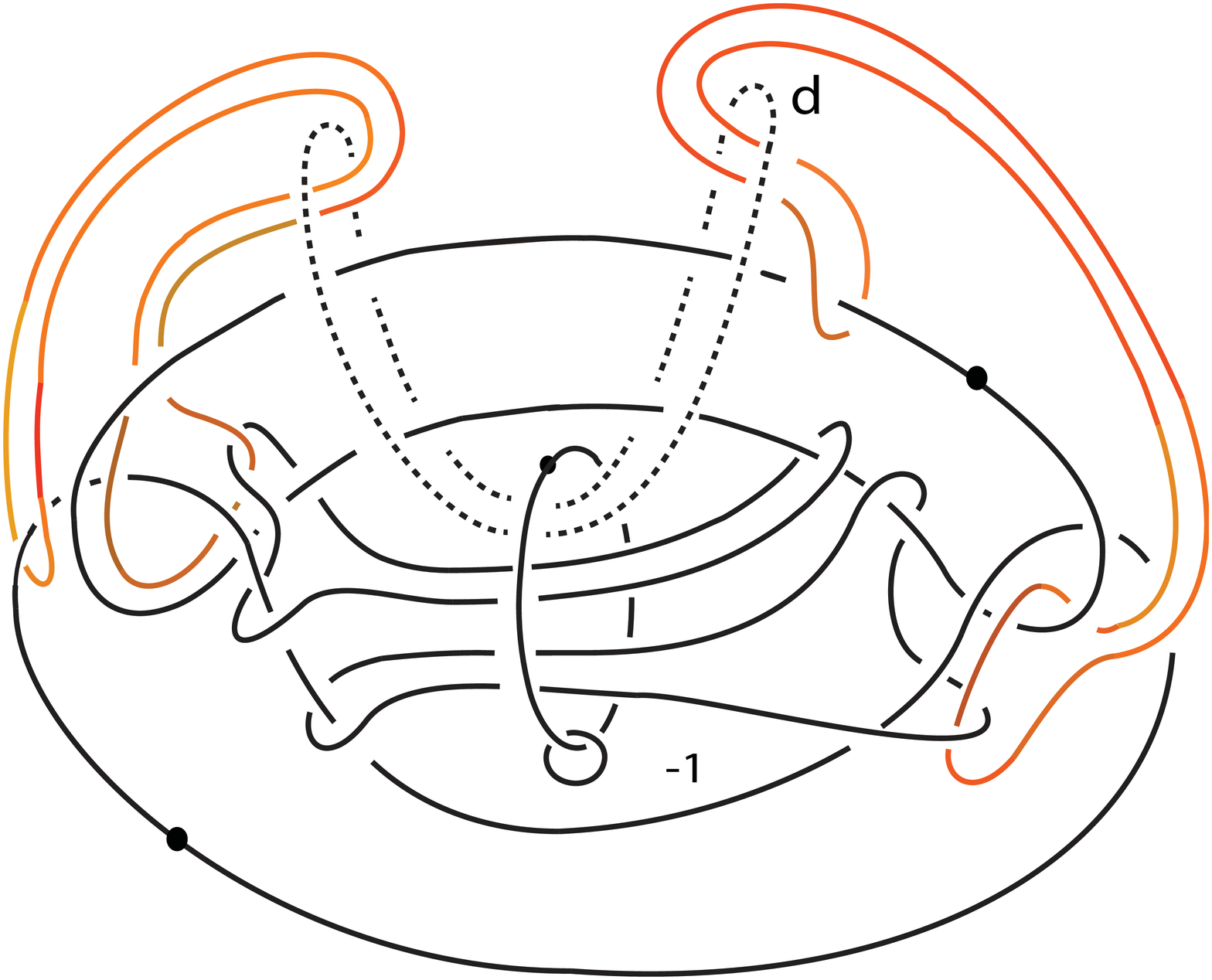}       
\caption{$d\subset \Sigma_{0}=B^4$}   \label{a6} 
\end{center}
 \end{figure}

\begin{figure}[ht]  \begin{center}
 \includegraphics[width=1\textwidth]{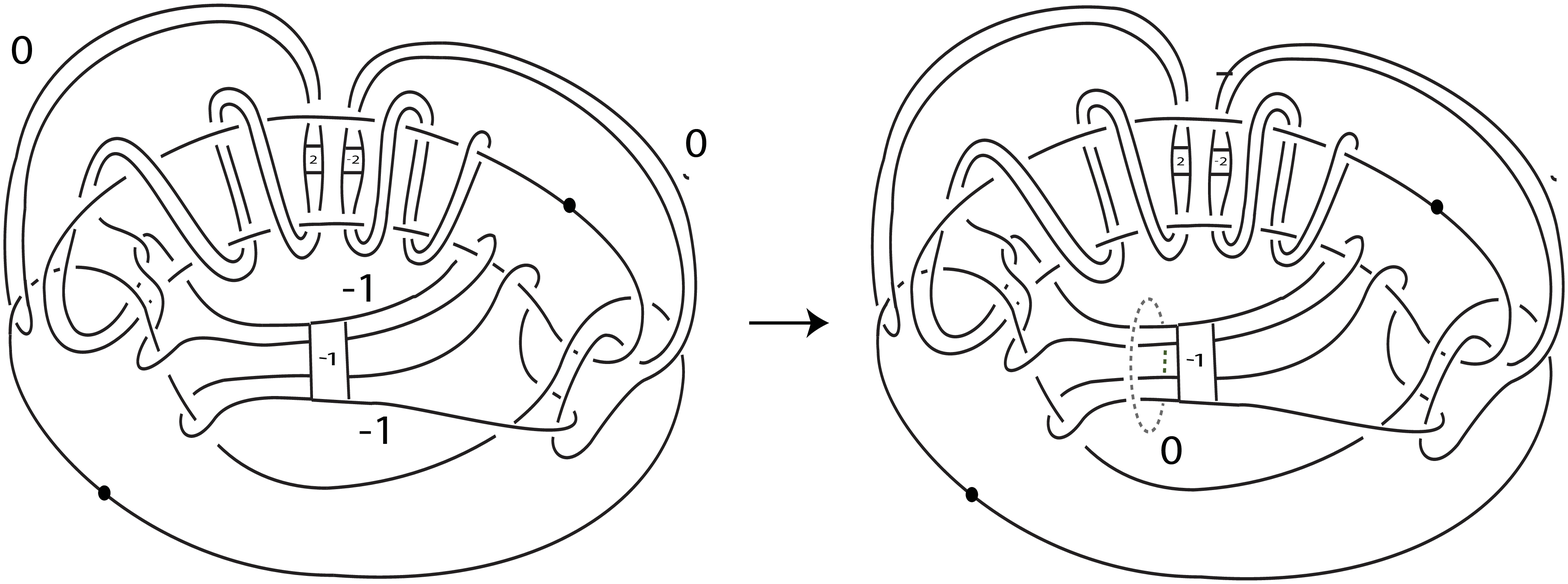}       
\caption{ $\Sigma_{n}$}     \label{a7} 
\end{center}
 \end{figure}

\begin{figure}[ht]  \begin{center}
 \includegraphics[width=1\textwidth]{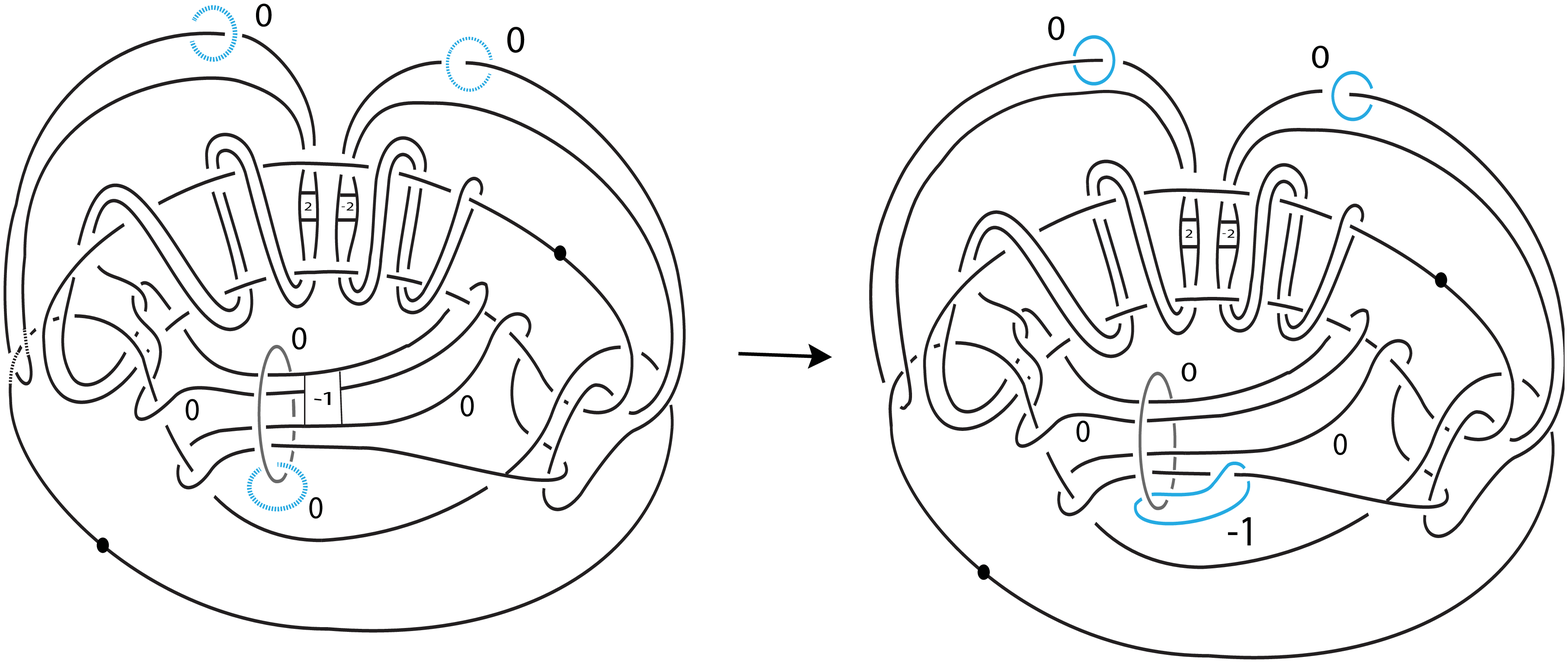}       
\caption{Turning $\Sigma_{n}$ upside down}      \label{a8} 
\end{center}
 \end{figure}

\begin{figure}[ht]  \begin{center}
 \includegraphics[width=1.1\textwidth]{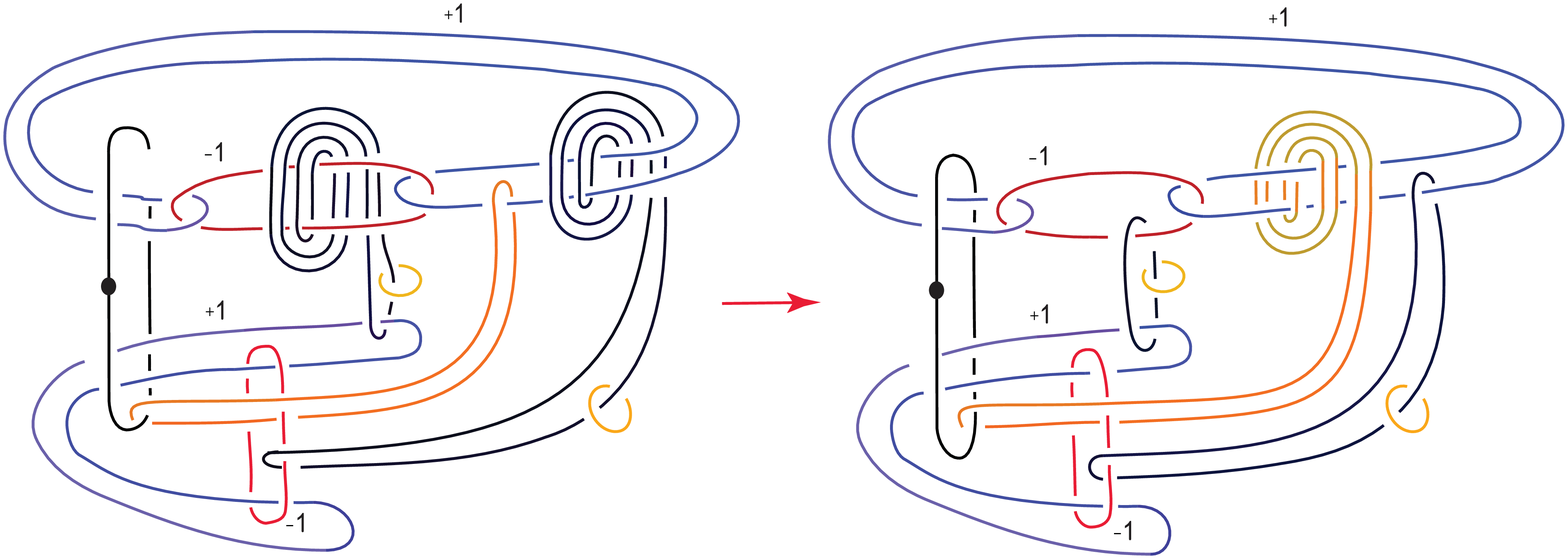}       
\caption{$\Sigma_{n}$, for n=2}     \label{a9} 
\end{center}
 \end{figure}

\begin{figure}[ht]  \begin{center}
 \includegraphics[width=.7\textwidth]{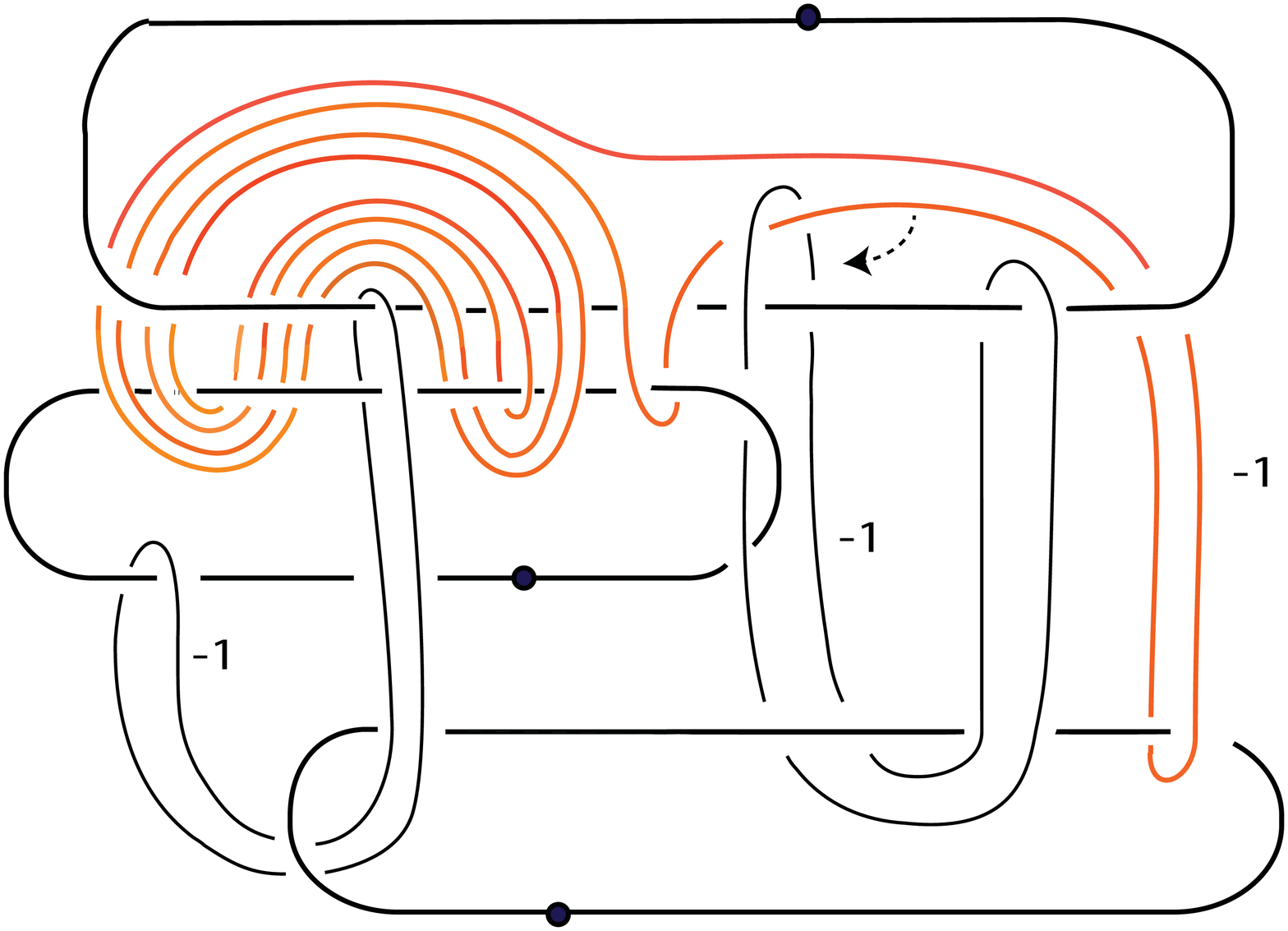}       
\caption{$\Sigma_{n}$, for $n=2$}      \label{a10} 
\end{center}
 \end{figure}

\begin{figure}[ht]  \begin{center}
 \includegraphics[width=.67\textwidth]{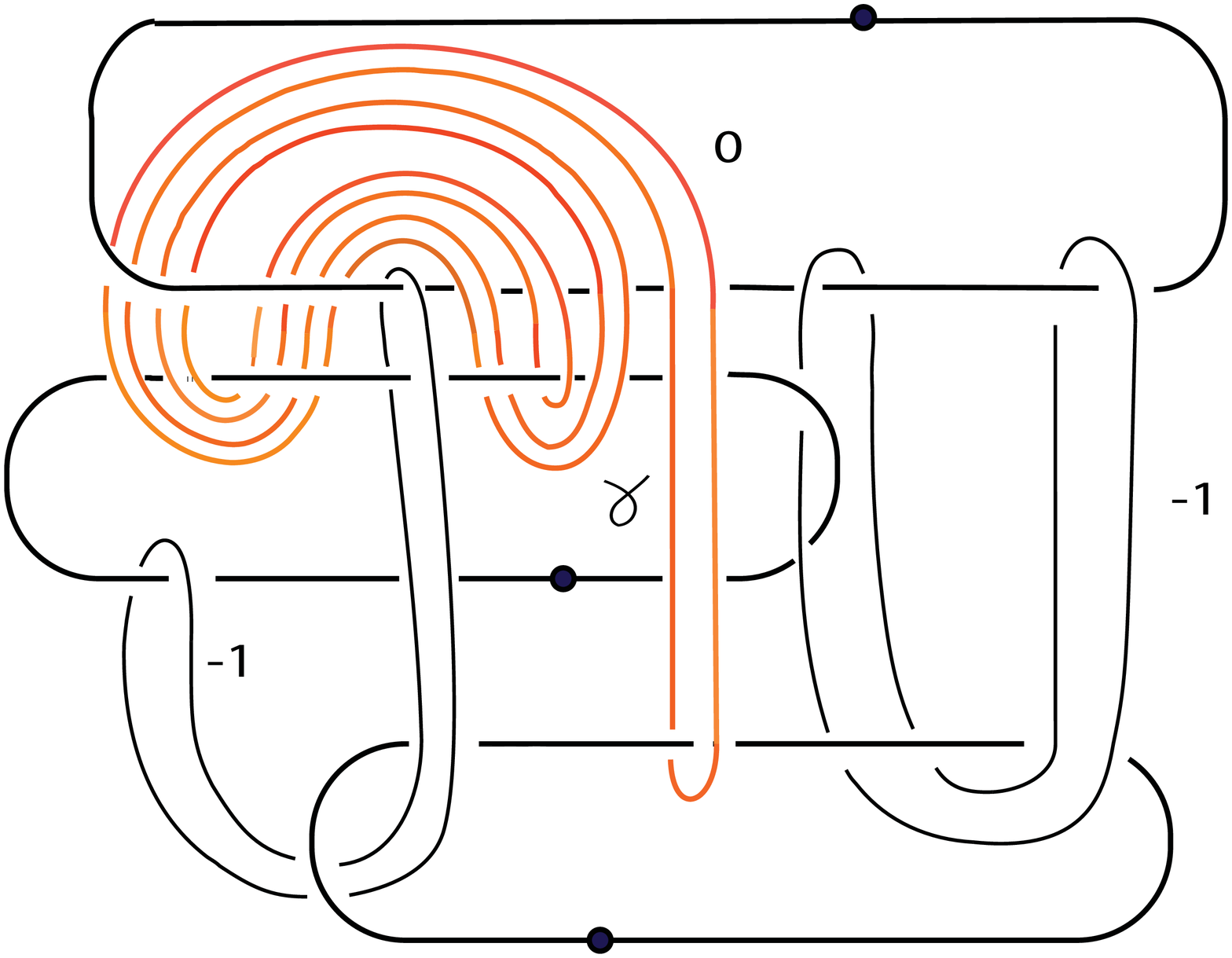}       
\caption{$\Sigma_{n}$, for $n=2$}      \label{a11} 
\end{center}
 \end{figure}

\begin{figure}[ht]  \begin{center}
 \includegraphics[width=.67\textwidth]{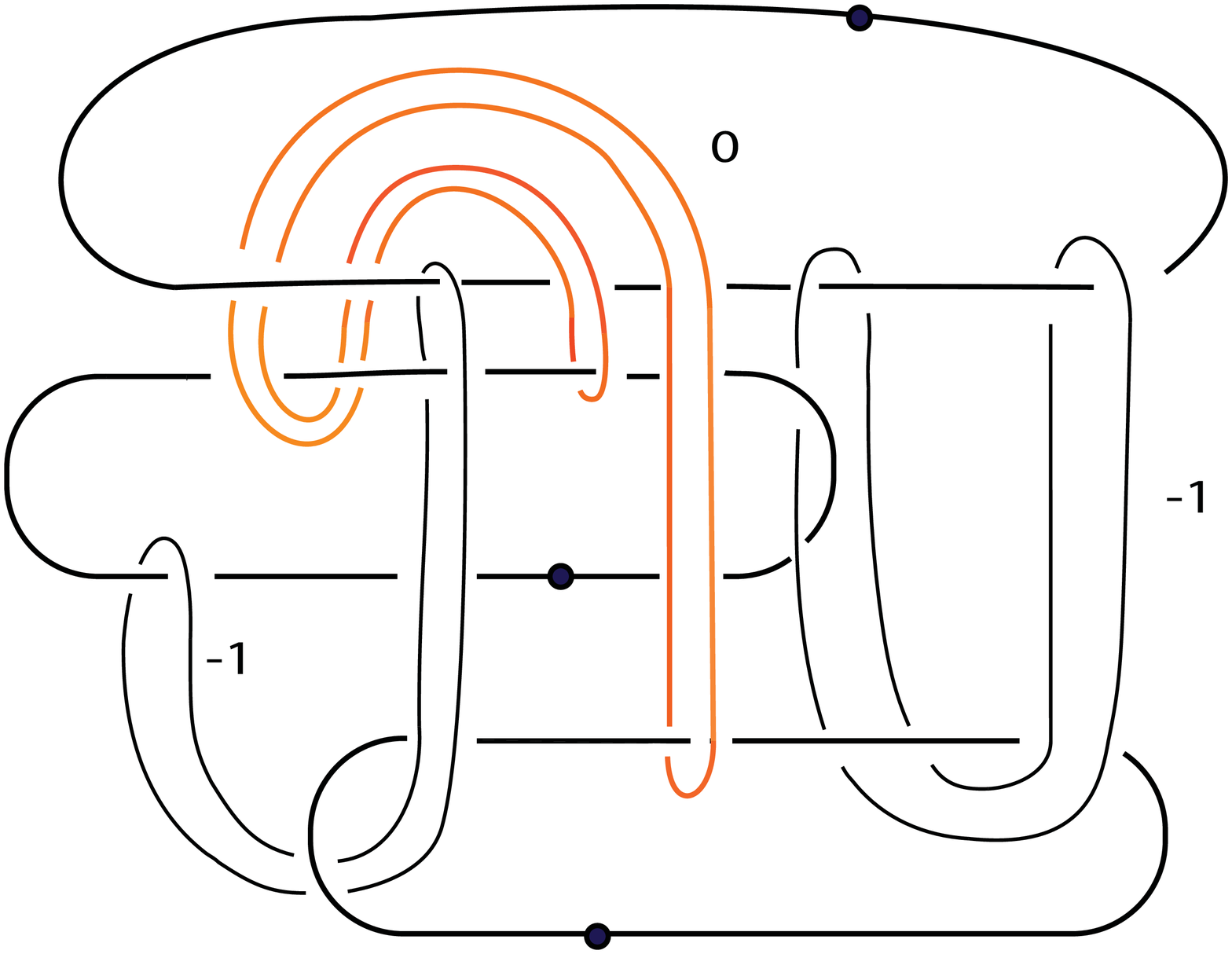}       
\caption{$\Sigma_{1}$}      \label{a12} 
\end{center}
 \end{figure}

\begin{figure}[ht]  \begin{center}
 \includegraphics[width=.7\textwidth]{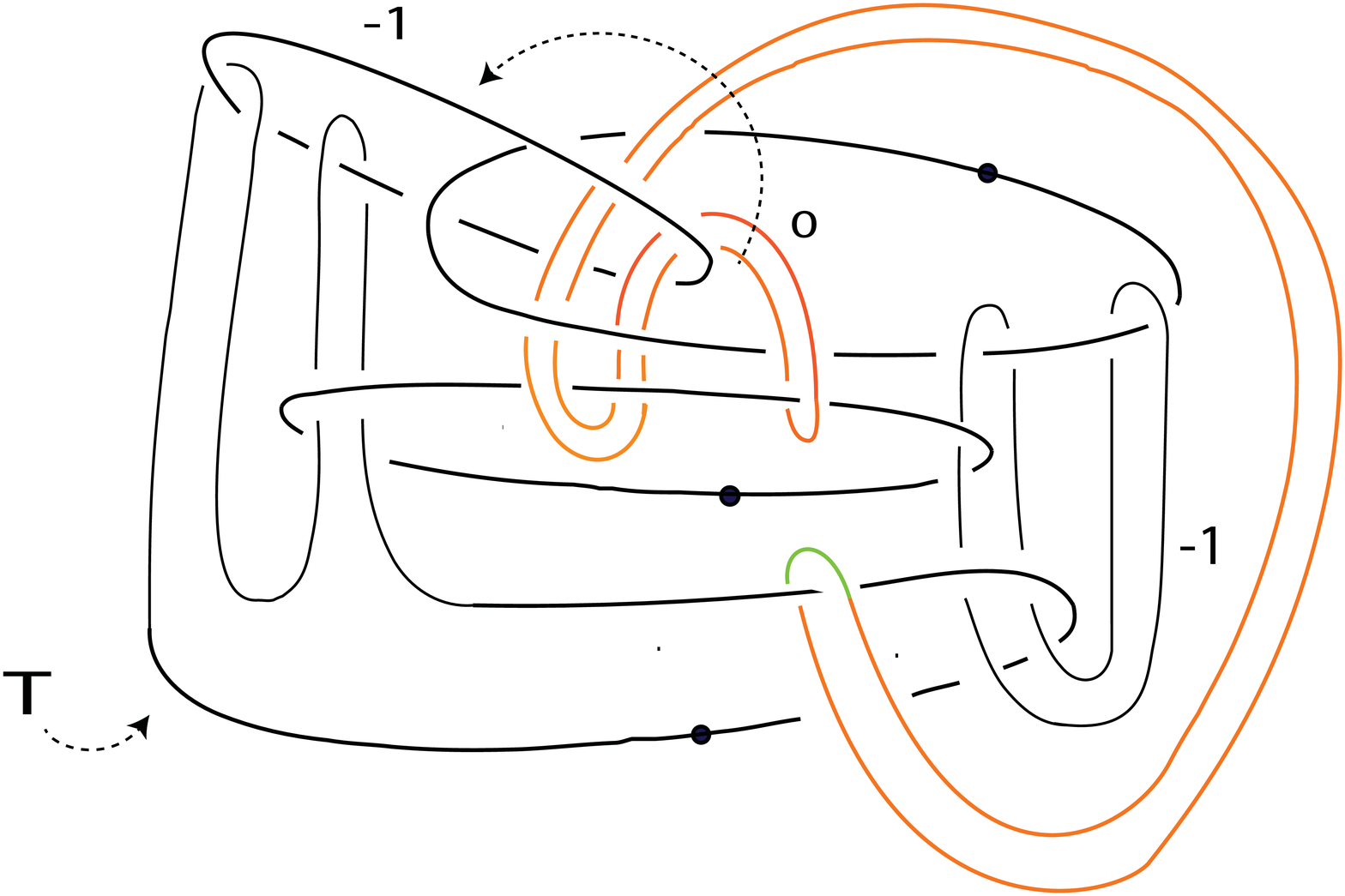}       
\caption{$\Sigma_{1}$}      \label{a13} 
\end{center}
 \end{figure}

We are now ready for the proof of Theorem~\ref{main}: Recall that we have the identification of the ribbon complements of Figure~\ref{a11} (without the curve $\gamma$), which is Figure~\ref{a19} (without the curve $\gamma$), with  Figure~\ref{a1} (without the small $-1$ linking circle). Now  a patient reader can easily check that under these identifications the curve $\gamma$ of Figure~\ref{a19} (which we also denote by 
$\gamma_{n}$) corresponds to the position of the curve curve $\gamma$  of Figure~\ref{a26}, after $\delta^{n} $ -move diffeomorphism (as described by Figure~\ref{a2}). It is amusing that $\delta$ of Figure~\ref{a26} remains unknot even after we attached a $2$-handle $\gamma$ to $W$, so in particular $\delta$-move commutes with attaching the $2$-handle handle $\gamma$  (this can be checked by using the Figure~\ref{a3} identification). Now remarks of the last paragraph of Section~\ref{introduction}  finishes the proof. \qed

  \begin{figure}[ht]  \begin{center}
 \includegraphics[width=.7\textwidth]{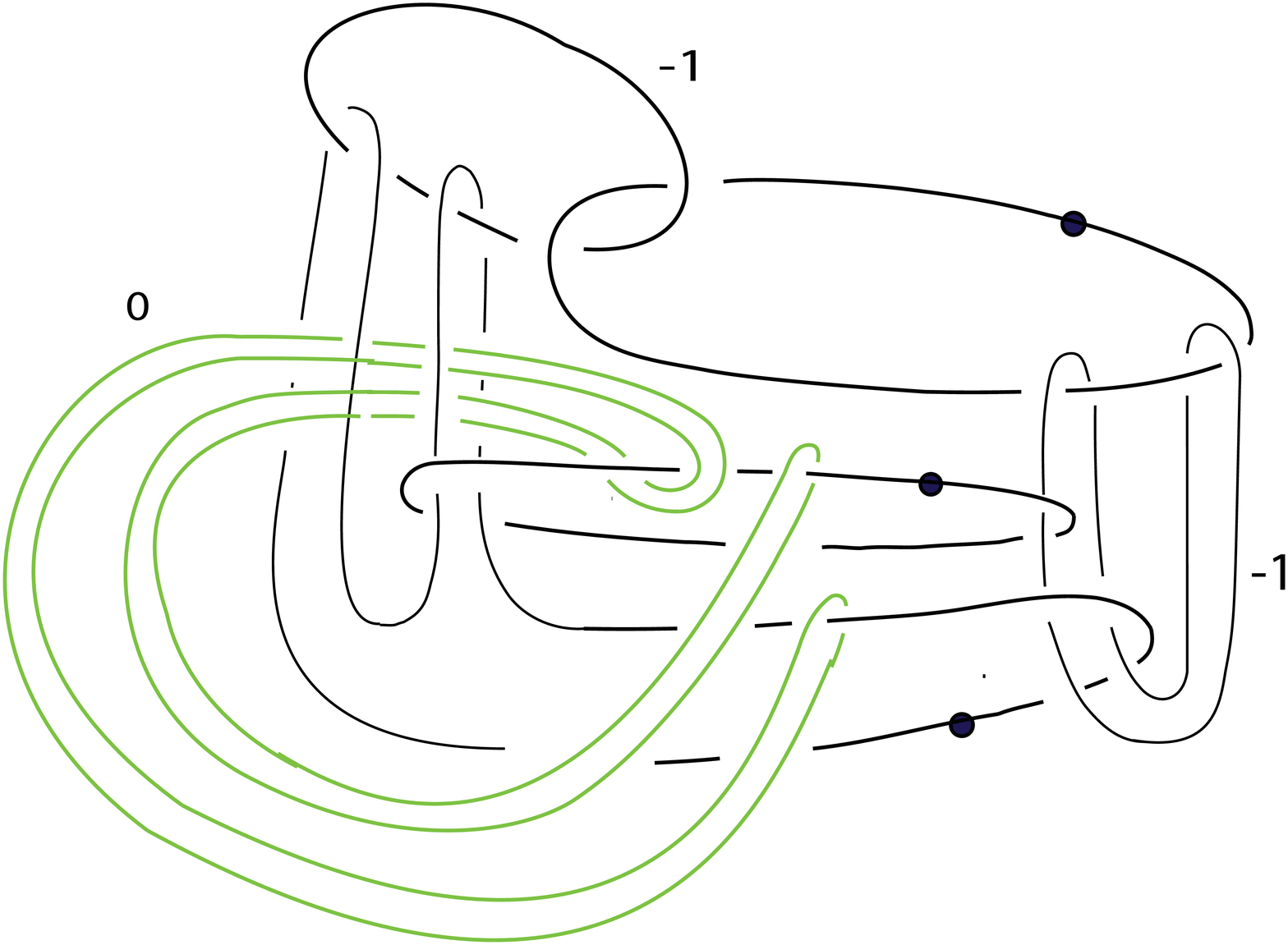}       
\caption{$\Sigma_{1}$}      \label{a14} 
\end{center}
 \end{figure}

 \begin{figure}[ht]  \begin{center}
 \includegraphics[width=.7\textwidth]{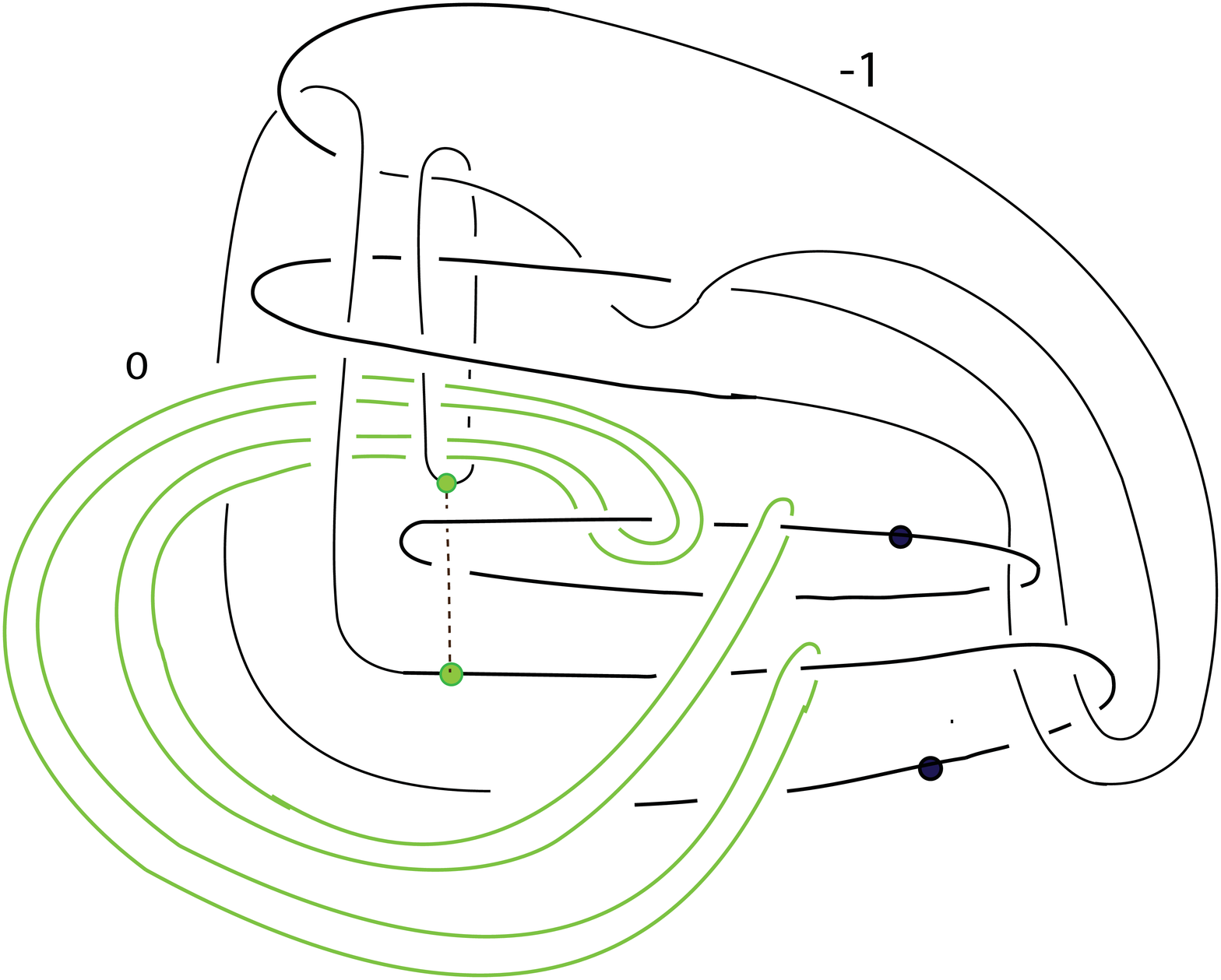}       
\caption{$\Sigma_{1}$}      \label{a15} 
\end{center}
 \end{figure}
 
  \begin{figure}[ht]  \begin{center}
 \includegraphics[width=.7\textwidth]{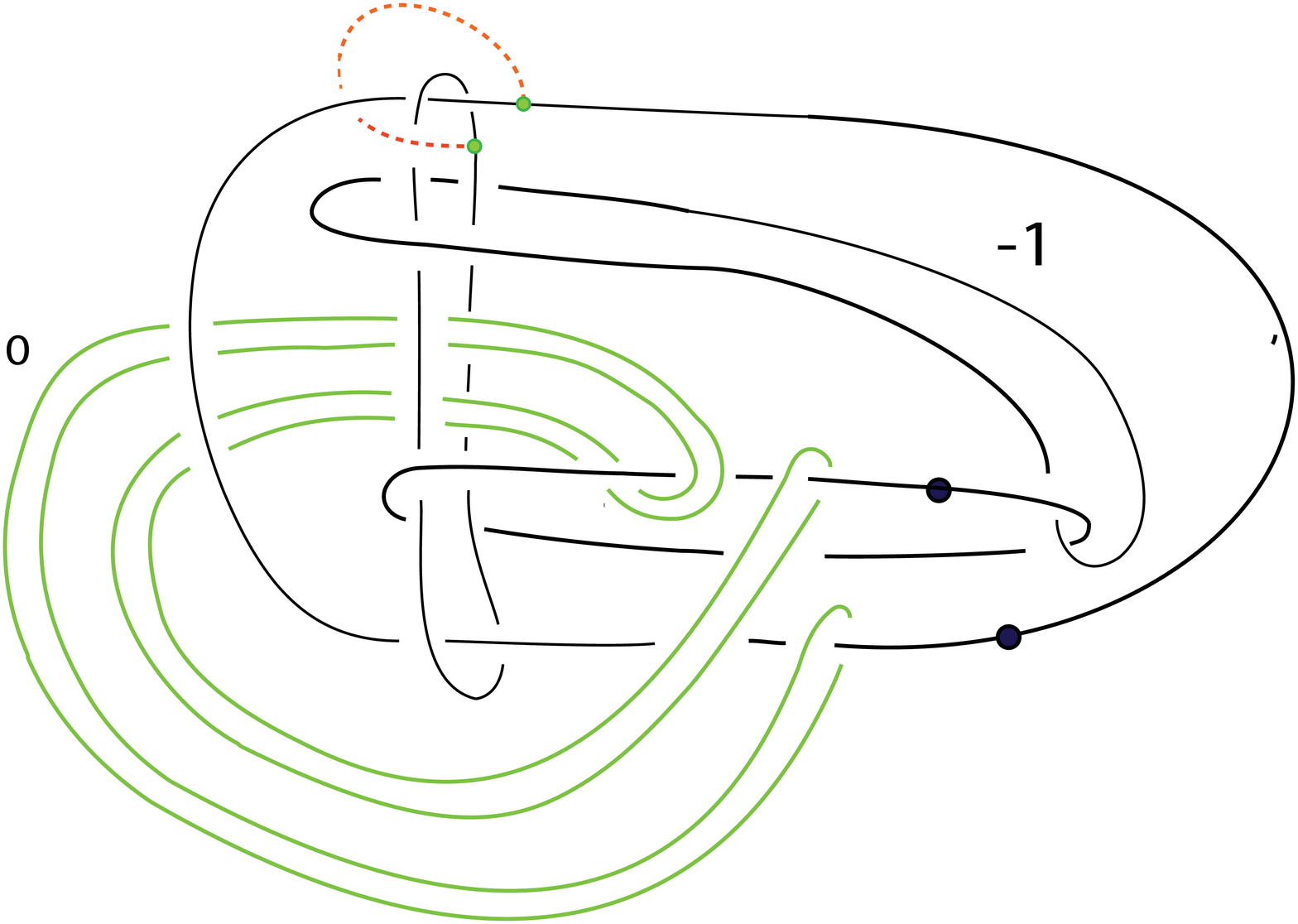}       
\caption{$\Sigma_{1}$}      \label{a16} 
\end{center}
 \end{figure}
 
  \begin{figure}[ht]  \begin{center}
 \includegraphics[width=.7\textwidth]{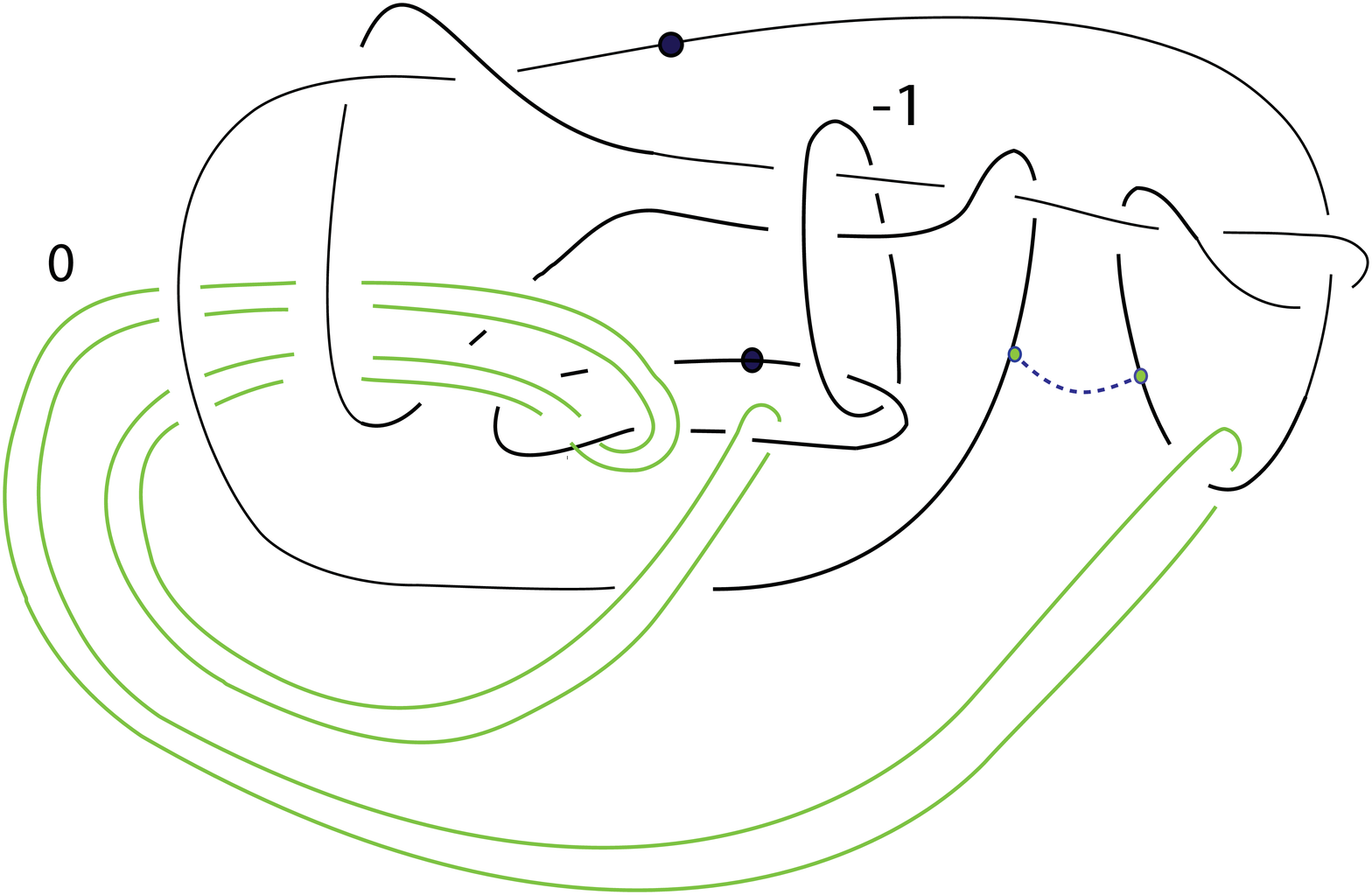}       
\caption{$\Sigma_{1}$}      \label{a17} 
\end{center}
 \end{figure}

 \begin{figure}[ht]  \begin{center}
 \includegraphics[width=.8\textwidth]{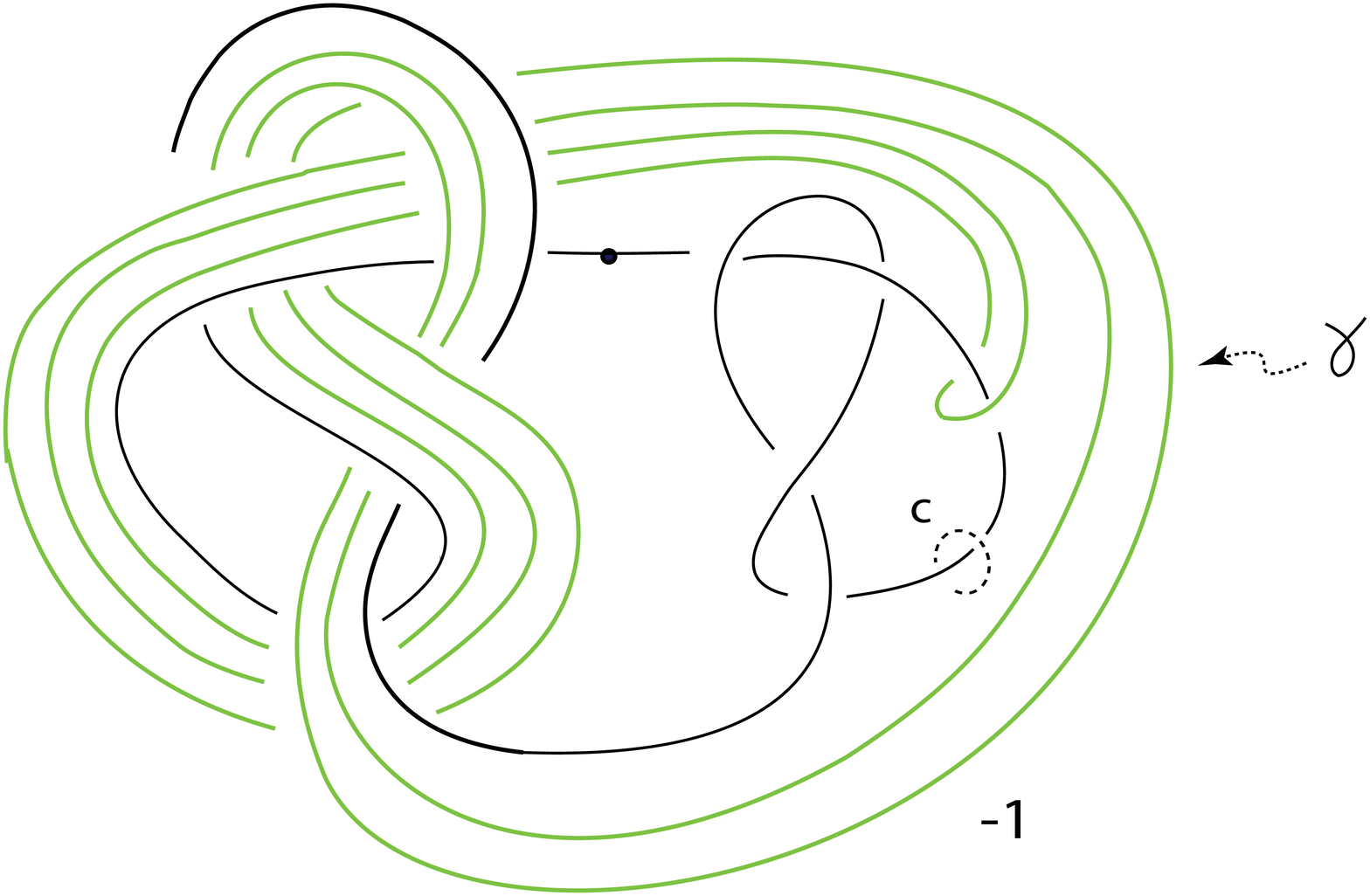}       
\caption{$\Sigma_{1}$}      \label{a18} 
\end{center}
 \end{figure}

 \begin{figure}[ht]  \begin{center}
 \includegraphics[width=.8\textwidth]{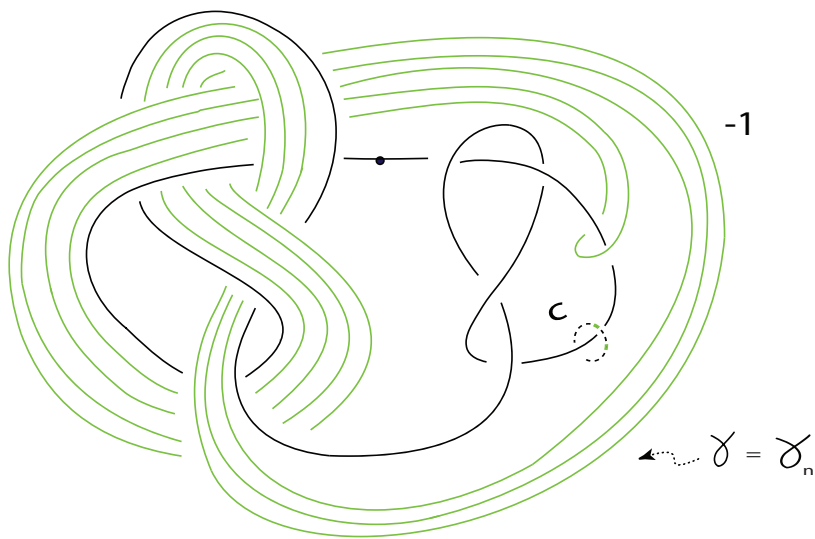}       
\caption{$\Sigma_{n}$, $n=2$}      \label{a19} 
\end{center}
 \end{figure}
 
 \begin{figure}[ht]  \begin{center}
 \includegraphics[width=.5\textwidth]{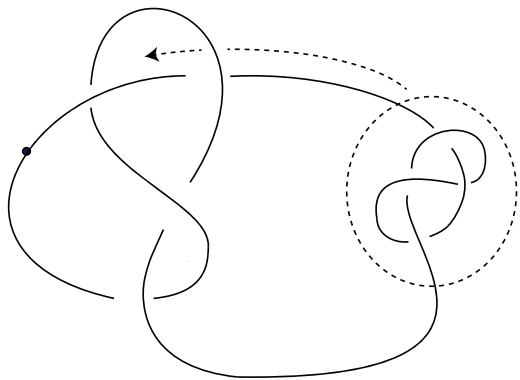}       
\caption{Rolling $f^{n}$}      \label{a20} 
\end{center}
 \end{figure}
 
 \begin{figure}[ht]  \begin{center}
 \includegraphics[width=.6\textwidth]{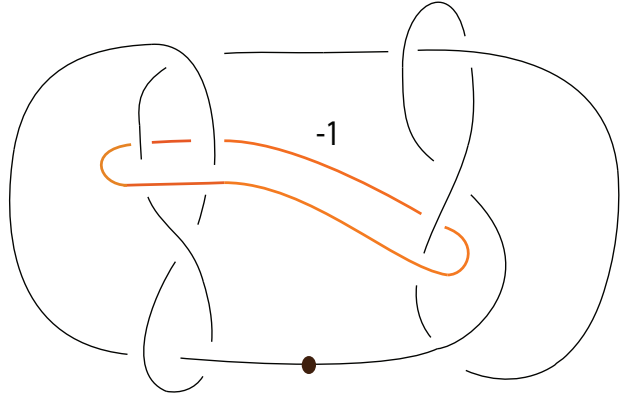}       
\caption{$\Sigma_{0}=B^4$}      \label{a21} 
\end{center}
 \end{figure}

 \begin{figure}[ht]  \begin{center}
 \includegraphics[width=1\textwidth]{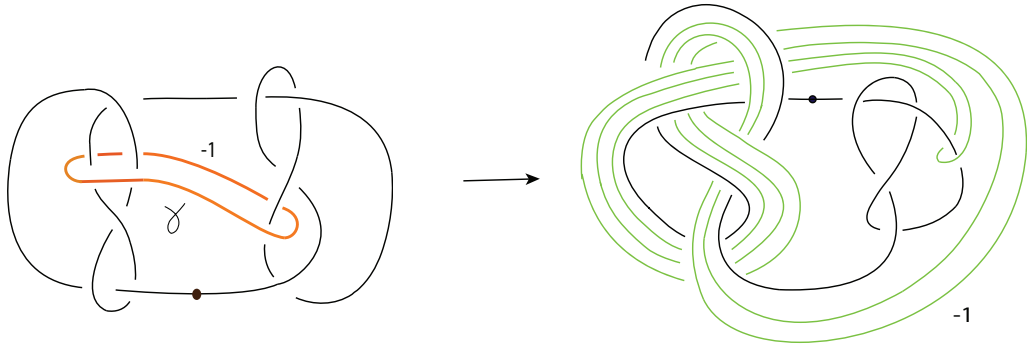}       
\caption{Rolling $2$-handle $\gamma$ by $f^{n}$}      \label{a22} 
\end{center}
 \end{figure}

\clearpage

\begin{figure}[ht]  \begin{center}
 \includegraphics[width=.85\textwidth]{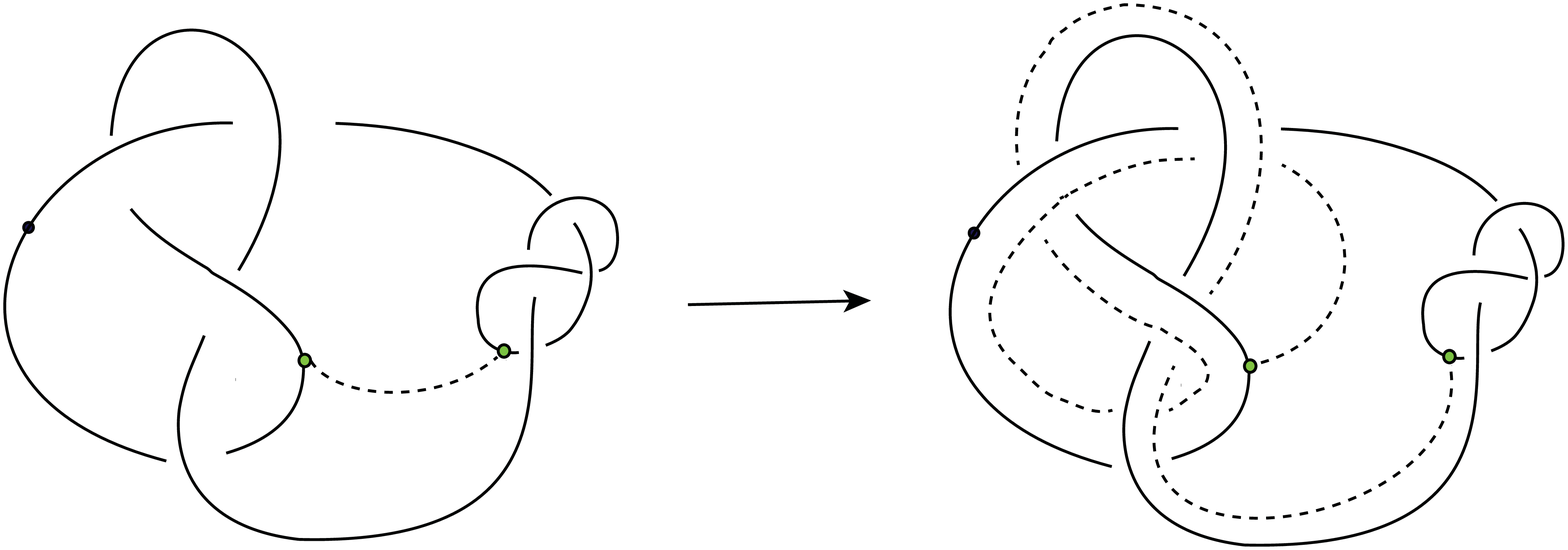}       
\caption{Carving ribbon $1$-handle by $f^{n}$}      \label{a23} 
\end{center}
 \end{figure}
 
 \begin{figure}[ht]  \begin{center}
 \includegraphics[width=1.1\textwidth]{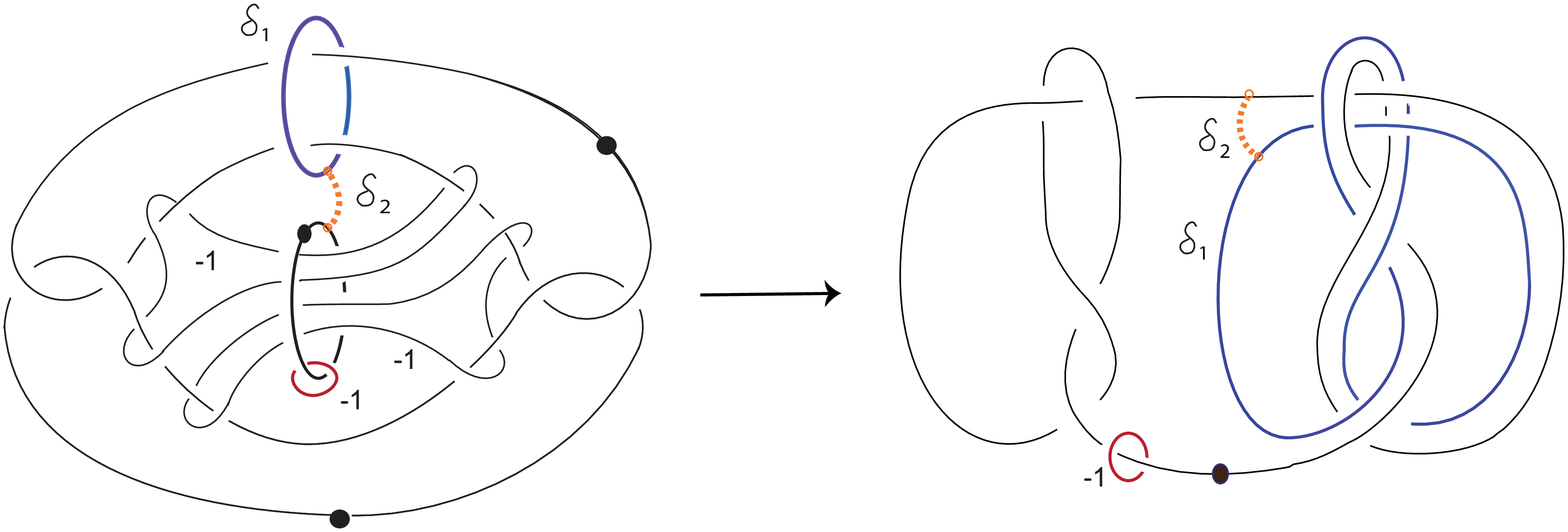}       
\caption{$\delta$-move $\rightsquigarrow $ Dehn surgery}      \label{a24} 
\end{center}
 \end{figure}

 \begin{figure}[ht]  \begin{center}
 \includegraphics[width=.55\textwidth]{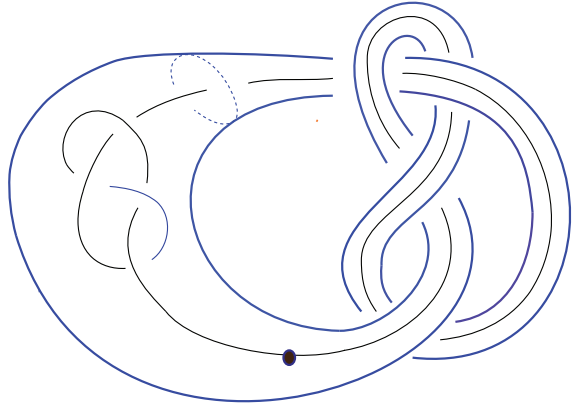}       
\caption{Dehn surgered torus}   \label{a25} 
\end{center}
 \end{figure}
 
   \begin{figure}[ht]  \begin{center}
 \includegraphics[width=1.1\textwidth]{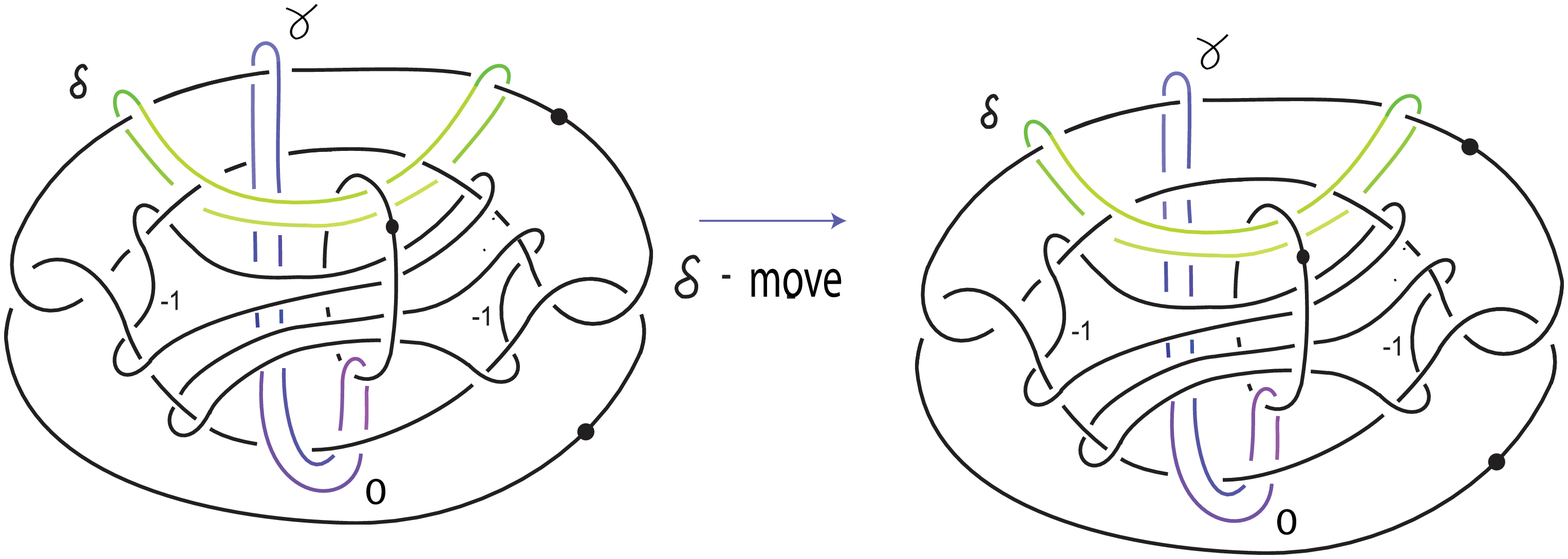}       
\caption{$\delta$-move}   \label{a26} 
\end{center}
 \end{figure}


\end{document}